\documentclass{article}

\usepackage{amssymb,latexsym,amsmath}

\usepackage{graphicx}

\usepackage{hyperref}

\begin{document}

\newcommand{\bfi}{\bfseries\itshape}

\makeatletter

\@addtoreset{figure}{section}

\def\thefigure{\thesection.\@arabic\c@figure}

\def\fps@figure{h, t}

\@addtoreset{table}{bsection}

\def\thetable{\thesection.\@arabic\c@table}

\def\fps@table{h, t}

\@addtoreset{equation}{section}

\def\theequation{\thesubsection.\arabic{equation}}

\makeatother

\newtheorem{theorem}{Theorem}[section]

\newtheorem{proposition}[theorem]{Proposition}

\newtheorem{lema}[theorem]{Lemma}

\newtheorem{cor}[theorem]{Corollary}

\newtheorem{definition}[theorem]{Definition}

\newtheorem{remark}[theorem]{Remark}

\newtheorem{exempl}{Example}[section]

\newenvironment{exemplu}{\begin{exempl}  \em}{\hfill $\square$

\end{exempl}}

\newcommand{\comment}[1]{\par\noindent{\raggedright\texttt{#1}

\par\marginpar{\textsc{Comment}}}}

\newcommand{\todo}[1]{\vspace{5 mm}\par \noindent \marginpar{\textsc{ToDo}}\framebox{\begin{minipage}[c]{0.95 \textwidth}

\tt #1 \end{minipage}}\vspace{5 mm}\par}

\newcommand{\ea}{\mbox{{\bf a}}}

\newcommand{\eu}{\mbox{{\bf u}}}

\newcommand{\ueu}{\underline{\eu}}

\newcommand{\ueo}{\overline{u}}

\newcommand{\oeu}{\overline{\eu}}

\newcommand{\ew}{\mbox{{\bf w}}}

\newcommand{\ef}{\mbox{{\bf f}}}

\newcommand{\eF}{\mbox{{\bf F}}}

\newcommand{\eC}{\mbox{{\bf C}}}

\newcommand{\en}{\mbox{{\bf n}}}

\newcommand{\eT}{\mbox{{\bf T}}}

\newcommand{\eL}{\mbox{{\bf L}}}

\newcommand{\eR}{\mbox{{\bf R}}}

\newcommand{\eV}{\mbox{{\bf V}}}

\newcommand{\eU}{\mbox{{\bf U}}}

\newcommand{\ev}{\mbox{{\bf v}}}

\newcommand{\eve}{\mbox{{\bf e}}}

\newcommand{\uev}{\underline{\ev}}

\newcommand{\eY}{\mbox{{\bf Y}}}

\newcommand{\eK}{\mbox{{\bf K}}}

\newcommand{\eP}{\mbox{{\bf P}}}

\newcommand{\eS}{\mbox{{\bf S}}}

\newcommand{\eJ}{\mbox{{\bf J}}}

\newcommand{\eB}{\mbox{{\bf B}}}

\newcommand{\eH}{\mbox{{\bf H}}}

\newcommand{\leb}{\mathcal{ L}^{n}}

\newcommand{\eI}{\mathcal{ I}}

\newcommand{\eE}{\mathcal{ E}}

\newcommand{\hen}{\mathcal{H}^{n-1}}

\newcommand{\eBV}{\mbox{{\bf BV}}}

\newcommand{\eA}{\mbox{{\bf A}}}

\newcommand{\eSBV}{\mbox{{\bf SBV}}}

\newcommand{\eBD}{\mbox{{\bf BD}}}

\newcommand{\eSBD}{\mbox{{\bf SBD}}}

\newcommand{\ecs}{\mbox{{\bf X}}}

\newcommand{\eg}{\mbox{{\bf g}}}

\newcommand{\paromega}{\partial \Omega}

\newcommand{\gau}{\Gamma_{u}}

\newcommand{\gaf}{\Gamma_{f}}

\newcommand{\sig}{{\bf \sigma}}

\newcommand{\gac}{\Gamma_{\mbox{{\bf c}}}}

\newcommand{\deu}{\dot{\eu}}

\newcommand{\dueu}{\underline{\deu}}

\newcommand{\dev}{\dot{\ev}}

\newcommand{\duev}{\underline{\dev}}

\newcommand{\weak}{\stackrel{w}{\approx}}

\newcommand{\mild}{\stackrel{m}{\approx}}

\newcommand{\lrightarrow}{\stackrel{L}{\rightarrow}}

\newcommand{\rrightarrow}{\stackrel{R}{\rightarrow}}

\newcommand{\strong}{\stackrel{s}{\approx}}

\newcommand{\weakdown}{\rightharpoondown}

\newcommand{\opg}{\stackrel{\mathfrak{g}}{\cdot}}

\newcommand{\opunu}{\stackrel{1}{\cdot}}
\newcommand{\opdoi}{\stackrel{2}{\cdot}}

\newcommand{\opn}{\stackrel{\mathfrak{n}}{\cdot}}
\newcommand{\opx}{\stackrel{x}{\cdot}}

\newcommand{\tr}{\ \mbox{tr}}

\newcommand{\Ad}{\ \mbox{Ad}}

\newcommand{\ad}{\ \mbox{ad}}

\renewcommand{\contentsname}{ }

\title{Boring mathematics,  \\
 {\it artistes pompiers} \\ 
 and impressionists}

\author{Marius Buliga \\
{\footnotesize Marius.Buliga@gmail.com}}

\date{This version:  06.09.2010}

\maketitle

In the middle of the 19th century, in France, just before the impressionist revolution,  painting was boring. Under the standards imposed by the Acad\'emie des Beaux Arts, paintings were done in a uniform technique, concerning a very restrictive list of subjects. 

I cite from 
\url{http://en.wikipedia.org/wiki/Impressionism#Beginnings}: 

\vspace{.5cm}

   "Colour was somber and conservative, and the traces of brush strokes were suppressed, concealing the artist's personality, emotions, and working techniques." 

\vspace{.5cm}

This seems very much similar with the situation in today's mathematical research. 

\vspace{.5cm}

Mind you, I am not claiming that mathematics is boring for the layman. This clich\'e is a cultural construct, 
but I am not interesting in it. I am referring to mathematical research! 

In any moment in history there are 
young men and women willing to do mathematical research, irrespective to the hostility of the cultural  
climate against it, or to the free fall of the quality of education. Nevertheless, from some time already, a good 
mathematician seems to be one with concealed personality, conservative (as in "follow the leader") goals 
and uniform working techniques. 

Here are some reasons why this is happening, in my opinion. 

\vspace{.5cm}

1. Because the mathematics  promoted these last decades IS BORING. To make a 
parallel with the  paintings universe, we have lots of Jean-L\'eon G\'er\^{o}me and Alexandre Cabanel, 
but very few C\'ezanne (well, Perelman fits the comparison with the painter) and for the moment 
we can't even dream about a Picasso!

\vspace{.5cm}

 2. If you want interesting mathematics then I politely suggest to look elsewhere, like in theoretical 
computer science, neuroscience or other emerging fields where simply breathtaking ideas, 
of mathematical flavor,  appear. And why not in mainstream mathematics? For the following two reasons. 

\vspace{.5cm}

 3. Too much politics! Since this mathematical research became a business, politics entered the stage. 
Even in this apparently rarefied world, research funds have to be attracted or managed, rankings have to 
be established, priorities have to be imposed. I am looking forward for an addition to Google Scholar, called maybe Scholar Mafia, offering for any given keyword a graph of clusters of researchers citing 
one another in the same cluster, by ignoring the other clusters... You know what I mean, moreover it is simple to be done. 

\vspace{.5cm}

4.  Brilliant is not enough! There could be an endless discussion if really a painter like Cabanel is 
worst than C\'ezanne, but it is somehow clear, now that the history filtered the events, that C\'ezanne had 
something genuinely new to show and Cabanel maybe not. Nevertheless, not enough brilliant masters 
will have lesser students and proteg\'es and this negative selection is in action from some time. 

\vspace{.5cm}

Please notice that the comparison between "established" mathematics and "l'art pompier" promoted by  19th century Acad\'emie des Beaux Arts is meaningful for yet another reason. The impressionist painters and the ones coming after them were not rejecting traditions. In fact they were more interested in the 
painting tradition than their more officially recognized colleagues. In mathematics, for example, geometry 
was butchered so much that for a while  mathematicians trained in "traditional geometry" could be found 
easier  in countries like Russia or Romania.  With the uniformization of the curriculum, this kind of 
mathematicians is in danger to dissapear. Nobody knows what our ancestor mathematician have done. 
Or, wait! this is not true, there is some information on the net, deformed by incompetence, but 
nevertheless available. 

This brings me to think about the future. Continuing the parallel with the impressionist revolution, 
maybe something like this will happen in mathematics as well, stimulated by the greatest invention 
in the recent history: the www. Indeed, maybe the ArXiv is to be compared with the "Salon des Refus\'es" 
which started  the painting revolution? Well, the vessel exists, it is now to be filled by really brilliant 
mathematicians.  Things should, and certainly will change regarding the way mathematics research is 
published and evaluated, too.

\vspace{.5cm}

Dear reader, please don't ask me who are the mathematicians corresponding to Cabanel, or other questions in the same vein. This note is not intended to have a key, or to polemically point to some 
person or institution. These are thoughts  that I hold since a long time and which I think they might be somehow relevant for explaining the situation of mathematical research,  in particular.

Briefly stated, I am looking forward for an impressionist revolution!

\end{document}